\documentclass{article}

\usepackage{geometry}
\usepackage{amsmath}
\usepackage{amssymb}
\usepackage{mathtools}
\usepackage{tikz}
\usepackage{hyperref}
\usepackage{pdfpages}
\hypersetup{colorlinks=true,linkcolor=blue}
\usepackage[hypcap=true]{caption}
\usepackage[thmmarks,thref,hyperref]{ntheorem}
\usepackage{thmtools}

\theoremstyle{plain}
\theoremseparator{. }
\newtheorem{theorem}{Theorem}
\newtheorem{lemma}[theorem]{Lemma}
\newtheorem{conjecture}[theorem]{Conjecture}
\newtheorem{claim}{Claim}

\def\titlevar{Circumference of essentially 4-connected planar triangulations}
\def\authorvar{Igor Fabrici, Jochen Harant, Samuel Mohr, Jens M. Schmidt}
\def\departA{Institute of Mathematics, P.J.~\v{S}af\'{a}rik University in Ko\v{s}ice, Slovakia}
\def\departB{Institute of Mathematics, Ilmenau University of Technology, Germany}
\def\departC{Institute for Algorithms and Complexity, Hamburg University of Technology, Germany}
\def\thedate{\today}

\title\titlevar
\author\authorvar
\date\thedate

\usepackage{hyperref}
\hypersetup{
    colorlinks=true,
    linkcolor=blue,
    citecolor=blue,
    pdftitle=\titlevar,
    pdfauthor=\authorvar
    bookmarks=true,
    pdfpagemode=UseOutlines,
}

\makeatletter

\newcommand{\ml}{l\kern-0.55mm\char39\kern-0.3mm}
\newcommand{\eF}{F_{\rm e}}
\newcommand{\ef}{f_{\rm e}}
\DeclareMathOperator{\circum}{circ}

\makeatother

\begin{document}
\begin{center}
{\bf \Large \titlevar }\\[8mm]
{\bf Igor~Fabrici\textsuperscript{\textnormal{1,2,a}},
 Jochen~Harant\textsuperscript{\textnormal{1,b}},
 Samuel~Mohr\textsuperscript{\textnormal{1,b}},
 Jens~M.~Schmidt\textsuperscript{\textnormal{1,c}}}\\[3mm]
\textsuperscript{a} {\it \departA }\\
\textsuperscript{b} {\it \departB }\\
\textsuperscript{c} {\it \departC }\\
\footnotetext[1]{Partially supported by DAAD, Germany (as part of BMBF) and the Ministry of Education, Science, Research and Sport of the Slovak
 Republic within the project 57447800.}
\footnotetext[2]{Partially supported by the Slovak Research and Development Agency under contract No.\ APVV-15-0116.}
\end{center}

\hrulefill

\begin{abstract}
A $3$-connected  graph $G$ is essentially $4$-connected if, for any $3$-cut $S\subseteq V(G)$ of $G$, at most one component of $G-S$ contains at
least two vertices. We prove that every essentially $4$-connected maximal planar graph $G$ on $n$ vertices contains a cycle of length at least
$\frac{2}{3}(n+4)$; moreover, this bound is sharp.
\end{abstract}


\smallskip
{\it $2010$ MSC:} 05C38, 05C10\\

\hrulefill

\vspace{7mm}


We consider finite, simple, and undirected graphs. The \textit{circumference} $\circum(G)$ of a graph $G$ is the length of a longest cycle of $G$. 
A cycle $C$ of  $G$ is an \textit{outer independent cycle} of $G$ if the set $V(G)\setminus V(C)$ is independent. (Note that an outer independent 
cycle is sometimes called a dominating cycle (\cite{Broersma2002}), although this is in contrast to the more commonly used definition of a dominating 
subgraph $H$ of $G$, where $V(H)$ dominates $V(G)$ in the usual sense.) A~set $S\subseteq V(G)$ ($S\subseteq E(G)$) is a \textit{$k$-cut} 
(a \textit{$k$-edge-cut}) of $G$ if $|S|=k$ and $G-S$ is disconnected. A~3-cut (a~3-edge-cut) $S$ of a 3-connected (3-edge-connected) graph $G$ is 
\textit{trivial} if at most one component of $G-S$ contains at least two vertices and the graph $G$ is \textit{essentially $4$-connected} 
(\textit{essentially $4$-edge-connected}) if every 3-cut (3-edge-cut) of $G$ is trivial. A 3-edge-connected graph $G$ is \textit{cyclically 
$4$-edge-connected} if for every 3-edge-cut $S$ of $G$, at most one component of $G-S$ contains a cycle.

It is well-known that for (3-connected) cubic graphs different from the triangular prism $K_3\times K_2$ (which is essentially 4-connected only) these 
three notions coincide (see e.g. \cite{Fleischner1989} and \cite{VanCleemput2018}). Obviously, the line graph $H=L(G)$ of a $3$-connected graph $G$ is 
4-connected if and only if $G$ is essentially 4-edge-connected. These two observations are  reasons for the quite great interest in studying all these 
three concepts of connectedness of graphs intensively.

Zhan \cite{Zhan1986} proved that every 4-edge-connected graph has a Hamiltonian line graph. Broersma \cite{Broersma2002} conjectured that even every 
essentially 4-edge-connected graph has a Hamiltonian line graph and showed that this is equivalent to the conjecture of Thomassen \cite{Thomassen1986} 
stating that every 4-connected line graph is Hamiltonian (which is known to be equivalent to the conjecture by Matthews and Sumner \cite{Matthews1985} 
stating that every 4-connected claw-free graph is Hamiltonian, as shown by Ryj\'{a}\v{c}ek \cite{Ryjacek1997}). Among others, the subclass of 
essentially 4-edge-connected cubic graphs is interesting due to a conjecture of Fleischner and Jackson \cite{Fleischner1989} stating that every 
essentially 4-edge-connected cubic graph has an outer independent cycle which is equivalent to the previous three conjectures.

Regarding to the existence of long cycles in essentially $4$-connected graphs we mention the following

\begin{conjecture}[Bondy, see \cite{Jackson1986}] \label{cBondy}
There exists a constant $c$, $0<c<1$, such that for every essentially $4$-connected cubic graph on $n$ vertices, $\circum(G)\ge cn$.
\end{conjecture}

Note that the conjecture of Fleischner and Jackson implies Conjecture~\ref{cBondy} with $c=\frac34$. Bondy's conjecture was later extended to all 
cyclically 4-edge-connected graphs (see \cite{Fleischner1989}). M\'{a}\v{c}ajov\'{a} and Maz\'{a}k \cite{Macajova2016} constructed essentially 
$4$-connected cubic graphs on $n=8m$ vertices with circumference $7m+2$. We remark that the conjecture of Fleischner and Jackson and, therefore, also 
Bondy's Conjecture  with $c=\frac34$ (this is the result of Gr\"unbaum and Malkevitch \cite{Grunbaum1976}) are true for planar graphs, which can be 
seen easily by the forthcoming Lemma~\ref{OI3C}. Many results concerning the circumference of  essentially $4$-connected planar graphs $G$ can be found 
in the literature.

For the class of essentially 4-connected cubic planar graphs, Tutte \cite{Tutte1960} showed that it contains a non-Hamiltonian graph, Aldred, Bau, 
Holton, and McKay \cite{Aldred2000} found a smallest non-Hamiltonian graph on 42 vertices, and Van Cleemput and Zamfirescu \cite{VanCleemput2018} 
constructed a non-Hamiltonian graph on $n$ vertices for all even $n\ge 42$. As already mentioned, Gr\"unbaum and Malkevitch \cite{Grunbaum1976} proved 
that $\circum(G)\ge \frac34n$ for any essentially 4-connected cubic planar graph $G$ on $n$ vertices and Zhang \cite{Zhang1987} (using the theory of 
Tutte paths) improved this lower bound on the circumference by 1. Recently, in \cite{Lo2018}, an infinite family of essentially 4-connected cubic 
planar graphs on $n$ vertices with circumference $\frac{359}{366}n$ was constructed.

In \cite{Jackson1992}, Jackson and Wormald extended the problem to find lower bounds on the circumference to the class of arbitrary essentially 
$4$-connected planar graphs. Their result $\circum(G)\ge\frac{2n+4}{5}$ was improved in \cite{Fabrici2020} to $\circum(G)\ge\frac58(n+2)$ for every 
essentially 4-connected planar graph $G$ on $n$ vertices. On the other side, there are infinitely many essentially 4-connected maximal planar graphs 
$G$ with $\circum(G)=\frac23(n+4)$ (\cite{Jackson1992}). To see this, let $G'$ be a 4-connected maximal planar graph on $n'\ge 6$ vertices and let $G$ 
be obtained from $G'$ by inserting a new vertex into each face of $G'$ and connecting it with all three boundary vertices of that face. Then $G$ is an 
essentially 4-connected maximal planar graph on $n=3n'-4$ vertices and, since $G'$ is Hamiltonian, it is easy to see that $\circum(G)=2n'=\frac23(n+4)$. 
It is still open whether there is an essentially 4-connected planar graph $G$ that satisfies $\circum(G)<\frac23(n+4)$. Indeed, we pose the following 
(to our knowledge so far unstated) Conjecture~\ref{con}, which has been the driving force in that area for over a decade.

\begin{conjecture}\label{con}
For every essentially $4$-connected planar graph on $n$ vertices, $\circum(G)\ge\frac23(n+4)$.
\end{conjecture}

By the forthcoming Theorem~\ref{main}, Conjecture~\ref{con} is shown to be true for  essentially 4-connected maximal planar graphs.

We remark that $G-S$ has exactly two components for every $3$-connected planar graph $G$ and every $3$-cut $S$ of $G$. Thus, in this case, $G$ is 
essentially $4$-connected if and only if $S$ forms the neighborhood of a vertex of degree $3$ of $G$ for every $3$-cut $S$ of $G$. This property will 
be used frequently in the proof of Theorem~\ref{main}.

A cycle $C$ of $G$ is a \textit{good cycle} of $G$ if $C$ is outer independent and $\deg_G(x)=3$ for all $x\in V(G)\setminus V(C)$. An edge $xy$ of 
a good cycle $C$ is \textit{extendable} if $x$ and $y$ have a common neighbor $z\in V(G)\setminus V(C)$. In this case, the cycle $C'$ of $G$, obtained 
from $C$ by replacing the edge $xy$ with the path $(x,z,y)$ is again good (and longer than $C$). The forthcoming Lemma~\ref{OI3C}  is an essential tool 
in the proof of Theorem~\ref{main} (an implicit proof for cubic essentially $4$-connected planar graphs can be found in \cite{Grunbaum1976}, the 
general case is proved in \cite{Fabrici2016}).

\begin{lemma} \label{OI3C}
Every essentially $4$-connected planar graph on $n\ge 11$ vertices contains a good cycle.
\end{lemma}

\begin{theorem} \label{main}
For every essentially $4$-connected maximal planar graph $G$ on $n\ge 8$ vertices,
$$\circum(G)\ge\frac{2}{3}(n+4).$$
\end{theorem}

\bigskip

\noindent
\textbf{Proof of Theorem~\ref{main}.}

Suppose $n\ge 11$, as for $n\in\{8,9,10\}$, Theorem~\ref{main} follows from the fact that $G$ is Hamiltonian (\cite{Barnette1970}). Using Lemma 
\ref{OI3C}, let $C=[v_1,v_2,\dots,v_k]$ (indices of vertices of $C$ are taken modulo $k$ in the whole paper) be a~longest good cycle of length $k$ of 
$G$ (i.e., $\circum(G)\ge k$) and let $H=G[V(C)]$ be the graph obtained from $G$ by removing all vertices of degree $3$ which do not belong to $C$. 
Obviously, $H$ is maximal planar and $C$ is a Hamiltonian cycle of $H$. A face $\varphi$ of $H$ is an \textit{empty face} of $H$ if $\varphi$ is also 
a face of $G$, otherwise $\varphi$ is a~\textit{non-empty face} of $H$. Denote by ${\eF}(H)$ the set of empty faces of $H$ and let $\ef(H)=|{\eF}(H)|$. 
Note that every face of $G$ has at least two (of three) vertices on $C$. The three neighbors of a vertex of $V(G)\setminus V(C)$ induce a separating 
3-cycle of $G$ creating the boundary of a non-empty face of $H$, which has no edge in common with $C$ because otherwise such an edge would be an 
extendable edge of $C$ in $G$.

Let $H_1$ and $H_2$ be the spanning subgraphs of $H$ consisting of the cycle $C$ and of its chords lying in the interior and in the exterior of $C$, 
respectively. Note that $E(H_1)\cap E(H_2)=E(C)$ and $H_1$ and $H_2$ are maximal outerplanar graphs, both having $k$-gonal outer face and $k-2$ 
triangular faces. Let $T_i$ be the weak dual of $H_i$, $i\in\{1,2\}$, which is the graph having all triangular faces of $H_i$ as vertex set such that 
two vertices of $T_i$ are adjacent if the triangular faces share an edge in $H_i$. Obviously, $T_i$ is a tree of maximum degree at most three.

A face $\varphi$ of $H$ is a \textit{$j$-face} if exactly $j$ of its three incident edges belong to $E(C)$. Since $n\ge 11$, there is no 3-face in $H$ 
and each face of $H$ is a $j$-face with $j\in\{0,1,2\}$. Denote by $f_j(H_i)$ the number of empty $j$-faces of $H_i$. Since $C$ does not contain any 
extendable edge, the following claim is obvious.

\begin{claim}\label{empty12}
Each face of $H$ incident with an edge of any longest good cycle (in particular, each $1$- or $2$-face) is empty.
\end{claim}

\noindent
An edge $e$ of $C$ incident with a $j$-face $\varphi$ and an $\ell$-face $\psi$, where $j,\ell\in\{1,2\}$, is a \textit{$(j,\ell)$-edge}. 
Let $\varphi$ be a 2-face of $H_i$. The sequence $B_{\varphi}=(\varphi,\varphi_2,\dots,\varphi_r)$, $r\ge 2$, is the \textit{$\varphi$-branch} if 
$\varphi_2,\dots,\varphi_{r-1}$ are 1-faces of $H_i$, $\varphi_{r}$ is a 0-face of $H_i$, and $\varphi_j$, $\varphi_{j+1}$ ($1\le j\le r-1$) are 
adjacent (i.e. $B_\varphi$ is a minimal path in $T_i$ with end vertices of degree 1 and 3). The \textit{rim} $R(B_\varphi)$ of the $\varphi$-branch 
$B_\varphi$ is the subgraph of $C$ induced by all edges of $C$ that are incident with an element of $B_\varphi$. Hence, it is easy to see:

\begin{claim}\label{rim path}
The rim of a $\varphi$-branch $B_{\varphi}=(\varphi,\varphi_2,\dots,\varphi_r)$ is a path of length $r$.
\end{claim}

\begin{claim} \label{empty0}
Let $\varphi=[v_1,v_2,v_3]$ be a $2$-face of $H_i$, let $B_\varphi=(\varphi,\varphi_2,\dots,\varphi_r)$, $r\ge 2$, be the $\varphi$-branch of $H_i$, 
and let $v_0v_2\in E(H_{3-i})$. If
\begin{enumerate}
\renewcommand{\labelenumi}{{\rm (\alph{enumi})}}
\item $R(B_\varphi)=(v_1,v_2,\dots,v_{r+1})$ is the rim of $B_\varphi$ or 
\item $R(B_\varphi)=(v_0,v_1,\dots,v_r)$ is the rim of $B_\varphi$ and $v_{-1}v_2\in E(H_{3-i})$, or 
\item $R(B_\varphi)=(v_{3-r},\dots,v_2,v_3)$ is the rim of $B_\varphi$ and $v_{-1}v_2\in E(H_{3-i})$, 
\end{enumerate}
then $\varphi_r$ is empty.
\end{claim}

\begin{center}
\begin{minipage}[t]{0.43\textwidth}
\begin{tikzpicture}[scale=0.07]
\draw[thick] (-15,0) -- (30,0);
\draw[thick] (40,0) -- (55,0);
\draw (-12,-4) node {$v_0$};
\draw (0,-4) node {$v_1$};
\draw (12,-4) node {$v_2$};
\draw (20,-4) node {$v_3$};
\draw (30,-4) node {$v_4$};
\draw (40,-4) node {$v_r$};
\draw (50,-4) node {$v_{r+1}$};
\draw (10,5) node {$\varphi$};
\draw (25,30) node {$\varphi_r$};
\draw[thick] (20cm,0cm) arc [start angle=0, end angle=180, radius=10cm];
\draw[thick] (30cm,0cm) arc [start angle=0, end angle=180, radius=15cm];
\draw[thick,cyan] (40cm,0cm) arc [start angle=0, end angle=180, radius=20cm];
\draw[thick,cyan] (50cm,0cm) arc [start angle=0, end angle=180, radius=25cm];
\draw[thick,cyan] (-10cm,0cm) arc [start angle=180, end angle=360, radius=10cm];
\draw[thick,cyan] (-15,0) -- (-10,0);
\draw[thick,cyan] (10,0) -- (30,0);
\draw[thick,cyan, dotted] (30,0) -- (40,0);
\draw[thick,cyan](50,0) -- (55,0);
\filldraw (-10,0) circle (8mm);
\filldraw (0,0) circle (8mm);
\filldraw (10,0) circle (8mm);
\filldraw (20,0) circle (8mm);
\filldraw (30,0) circle (8mm);
\filldraw (40,0) circle (8mm);
\filldraw (50,0) circle (8mm);
\filldraw[white] (20,-20) circle (5mm);
\end{tikzpicture}
\vspace{-8mm}
\begin{center}
(a)
\end{center}
\end{minipage}
\begin{minipage}[t]{0.52\textwidth}
\begin{tikzpicture}[scale=0.07]
\draw[thick] (-25,0) -- (30,0);
\draw[thick] (40,0) -- (55,0);
\draw (-24,-4) node {$v_{-1}$};
\draw (-12,-4) node {$v_0$};
\draw (0,-4) node {$v_1$};
\draw (12,-4) node {$v_2$};
\draw (20,-4) node {$v_3$};
\draw (30,-4) node {$v_4$};
\draw (40,-4) node {$v_{r-1}$};
\draw (50,-4) node {$v_r$};
\draw (10,5) node {$\varphi$};
\draw (20,35) node {$\varphi_r$};
\draw[thick] (20cm,0cm) arc [start angle=0, end angle=180, radius=10cm];
\draw[thick] (30cm,0cm) arc [start angle=0, end angle=180, radius=15cm];
\draw[thick] (50cm,0cm) arc [start angle=0, end angle=180, radius=25cm];
\draw[thick] (-10cm,0cm) arc [start angle=180, end angle=360, radius=10cm];
\draw[thick,cyan] (40cm,0cm) arc [start angle=0, end angle=180, radius=20cm];
\draw[thick,cyan] (50cm,0cm) arc [start angle=0, end angle=180, radius=30cm];
\draw[thick,cyan] (-20cm,0cm) arc [start angle=180, end angle=360, radius=15cm];
\draw[thick,cyan] (-25,0) -- (-20,0);
\draw[thick,cyan] (-10,0) -- (0,0);
\draw[thick,cyan] (10,0) -- (30,0);
\draw[thick,cyan, dotted] (30,0) -- (40,0);
\draw[thick,cyan] (50,0) -- (55,0);
\filldraw (-20,0) circle (8mm);
\filldraw (-10,0) circle (8mm);
\filldraw (0,0) circle (8mm);
\filldraw (10,0) circle (8mm);
\filldraw (20,0) circle (8mm);
\filldraw (30,0) circle (8mm);
\filldraw (40,0) circle (8mm);
\filldraw (50,0) circle (8mm);
\filldraw[white] (20,-20) circle (5mm);
\end{tikzpicture}
\vspace{-8mm}
\begin{center}
(b)
\end{center}
\end{minipage}\\[3mm]
\begin{minipage}[t]{0.7\textwidth}
\begin{tikzpicture}[scale=0.07]
\draw[thick] (-25,0) -- (30,0);
\draw[thick] (40,0) -- (60,0);
\draw[thick] (70,0) -- (85,0);
\draw (-24,-4) node {$v_{-1}$};
\draw (-12,-4) node {$v_0$};
\draw (0,-4) node {$v_1$};
\draw (12,-4) node {$v_2$};
\draw (20,-4) node {$v_3$};
\draw (30,-4) node {$v_4$};
\draw (40,-4) node {$v_{s-1}$};
\draw (50,-4) node {$v_s$};
\draw (60,-4) node {$v_{s+1}$};
\draw (70,-4) node {$v_{r-1}$};
\draw (80,-4) node {$v_r$};
\draw (10,5) node {$\varphi$};
\draw (35,50) node {$\varphi_r$};
\draw[thick] (30cm,0cm) arc [start angle=0, end angle=180, radius=15cm];
\draw[thick] (40cm,0cm) arc [start angle=0, end angle=180, radius=20cm];
\draw[thick] (50cm,0cm) arc [start angle=0, end angle=180, radius=25cm];
\draw[thick] (50cm,0cm) arc [start angle=0, end angle=180, radius=30cm];
\draw[thick] (60cm,0cm) arc [start angle=0, end angle=180, radius=35cm];
\draw[thick] (-10cm,0cm) arc [start angle=180, end angle=360, radius=10cm];
\draw[thick,cyan] (20cm,0cm) arc [start angle=0, end angle=180, radius=10cm];
\draw[thick,cyan] (70cm,0cm) arc [start angle=0, end angle=180, radius=40cm];
\draw[thick,cyan] (80cm,0cm) arc [start angle=0, end angle=180, radius=45cm];
\draw[thick,cyan] (-20cm,0cm) arc [start angle=180, end angle=360, radius=15cm];
\draw[thick,cyan] (-25,0) -- (-20,0);
\draw[thick,cyan] (0,0) -- (10,0);
\draw[thick,cyan] (20,0) -- (30,0);
\draw[thick,cyan, dotted] (30,0) -- (40,0);
\draw[thick,cyan] (40,0) -- (60,0);
\draw[thick,cyan, dotted] (60,0) -- (70,0);
\draw[thick,cyan] (80,0) -- (85,0);
\filldraw (-20,0) circle (8mm);
\filldraw (-10,0) circle (8mm);
\filldraw (0,0) circle (8mm);
\filldraw (10,0) circle (8mm);
\filldraw (20,0) circle (8mm);
\filldraw (30,0) circle (8mm);
\filldraw (40,0) circle (8mm);
\filldraw (50,0) circle (8mm);
\filldraw (60,0) circle (8mm);
\filldraw (70,0) circle (8mm);
\filldraw (80,0) circle (8mm);
\filldraw[white] (20,-17.32) circle (5mm);
\end{tikzpicture}
\vspace{-5mm}
\begin{center}
(c)
\end{center}
\end{minipage}
\bigskip

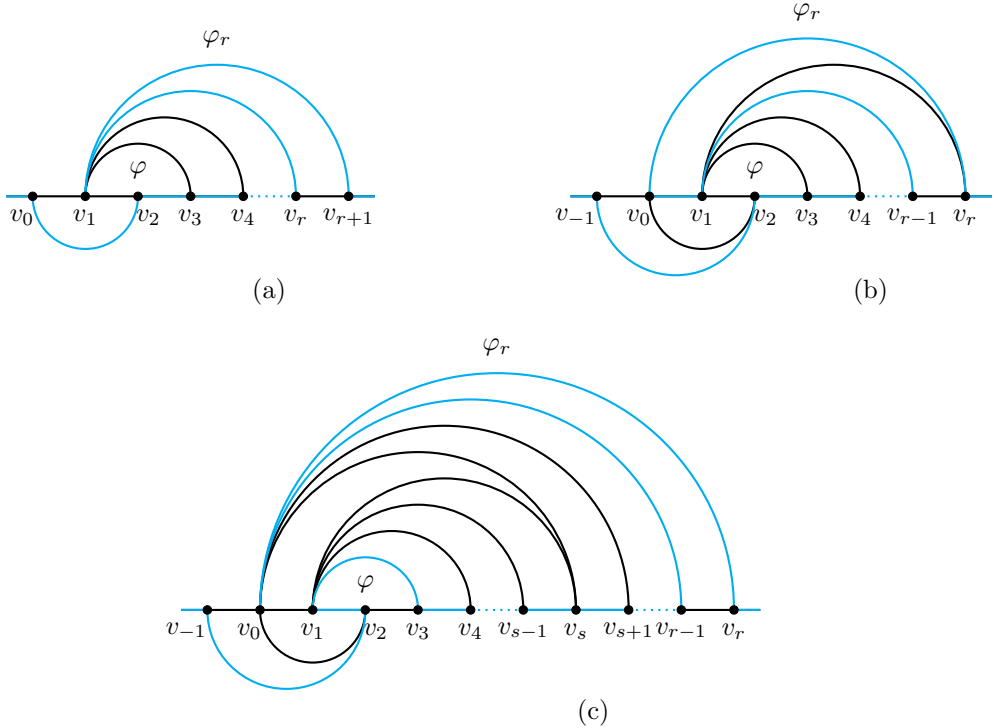
Fig.~1. A longest good cycle (cyan) sharing an edge with $\varphi_r$.
\end{center}

\newpage

\noindent
\textit{Proof.}

\textbf{(a)} The cycle $C'$ obtained from $C$ by replacing the path $(v_0,v_1,\dots,v_{r+1})$ with the path $(v_0,v_2,\dots,v_r,v_1,$ $v_{r+1})$ 
(Fig.~1(a)) is another longest good cycle of $G$ and contains the edge $v_1v_{r+1}$ incident with $\varphi_r$, thus $\varphi_r$ is empty (by 
Claim~\ref{empty12}).

\textbf{(b)} Let $\varphi_s=[v_0,v_1,v_{s}]$, for some $s$ with $3\le s\le r$, be a 1-face of $H_i$. The cycle $C'$ obtained from $C$ by replacing 
the path $(v_{-1},v_0,\dots,v_r)$ by the path $(v_{-1},v_2,\dots,v_{r-1},v_1,v_0,v_r)$, for $s=r$ (Fig.~1(b)), or by the path 
$(v_{-1},v_2,v_1,v_3,\dots,v_{r-1},v_0,v_r)$, for $s\le r-1$ (Fig.~1(c)), is a~longest good cycle of $G$ and contains the edge $v_0v_r$ 
incident with $\varphi_r$, thus $\varphi_r$ is empty (by Claim~\ref{empty12}).

\textbf{(c)} If $r\le 3$, then $\varphi_r$ is empty by (a) or (b). If $r\ge 4$, then $v_0v_3, v_{-1}v_3\in E(H_i)$, thus 
$\{v_{-1},v_2,v_3\}$ is a non-trivial 3-cut, a contradiction.
\hfill $\Box$

\bigskip
These tools will be used continuously in the following; we continue with the proof of Theorem~\ref{main}. Hereby, we consider two cases. In the first 
case, both subgraphs $H_1$ and $H_2$ have some 0-faces. By using a customized discharging method, we distribute some weights from edges to faces to 
prove that sufficiently many faces are empty (each empty face will finally contain weight at most $\frac23$). In the second case, there are only empty 
faces on one side of $C$, so that all vertices not in $C$ are located on the other side of $C$. We have to prove that there are some additional empty 
faces on this side.

\medskip
\noindent
\textbf{CASE 1.} Let $H_1$ and $H_2$ both contain at least two 0-faces or one non-empty 0-face.\\[2mm]
For every edge $e$ of $C$ we define the weight $w_0(e)=1$. Obviously, $\sum\limits_{e\in E(C)}w_0(e)=|E(C)|=k$.

\medskip
\noindent
\textbf{First redistribution of weights.}\\[2mm]
Each edge of $C$ sends weight to both incident faces as follows

\medskip
\noindent
\hspace*{-7pt}
{\setlength{\tabcolsep}{3pt}
\begin{tabular}{lp{12.2cm}}
\textbf{Rule R1.}& A (1,1)-edge sends $\frac12$ to both incident 1-faces.\\[3pt]
\textbf{Rule R2.}& A (1,2)-edge sends $\frac23$ to the incident 1-face and $\frac13$ to the incident 2-face.\\[3pt]
\textbf{Rule R3.}& A (2,2)-edge sends $\frac12$ to both incident 2-faces.
\end{tabular}}\\[3mm]

The edges of $C$ completely redistribute their weights to incident 1- and 2-faces. For an empty face $\varphi$, let $w_1(\varphi)$ be the total 
weight obtained by $\varphi$ (in first redistribution). Obviously, for an empty face $\varphi$, it is 
$$w_1(\varphi)=\left\{
\begin{array}{cl}
1,& \mbox{if $\varphi$ is a 2-face incident with two (2,2)-edges,}\\[3pt]
\frac56,& \mbox{if $\varphi$ is a 2-face incident with a (1,2)-edge and a (2,2)-edge,}\\[3pt]
\frac23,& \mbox{if $\varphi$ is a 2-face incident with two (1,2)-edges,}\\[3pt]
\frac23,& \mbox{if $\varphi$ is a 1-face incident with a (1,2)-edge,}\\[3pt]
\frac12,& \mbox{if $\varphi$ is a 1-face incident with a (1,1)-edge,}\\[3pt]
0,& \mbox{if $\varphi$ is a 0-face.}
\end{array}\right.
$$
Moreover, $\sum\limits_{\varphi\in\eF(H)} w_1(\varphi)=|E(C)|=k$.

\medskip
\newpage
\noindent
\textbf{Second redistribution of weights.}\\[2mm]
The weight of 2-faces of $H$ exceeding $\frac23$ will be redistributed to 1-faces and empty 0-faces of $H$ by the following rules. 
Let $\varphi$ be a 2-face of $H_i$ with $w_1(\varphi)>\frac23$ (i.e. incident with at least one (2,2)-edge) and let 
$B_{\varphi}=(\varphi,\varphi_2,\dots,\varphi_r)$, $r\ge 2$, be the $\varphi$-branch. Moreover, let $\alpha$ be a 2-face of $H_{3-i}$ 
adjacent to $\varphi$ and let $\alpha_2$ be the face of $H_{3-i}$ adjacent to $\alpha$.

\medskip
\noindent
\hspace*{-7pt}
{\setlength{\tabcolsep}{3pt}
\begin{tabular}{lp{14cm}}
\textbf{Rule R4.}& $\varphi$ sends $w_1(\varphi)-\frac23$ to $\varphi_r$ if $\varphi_r$ is empty and $r\le 3$.\\[3pt]
\textbf{Rule R5.}& $\varphi$ sends $\frac16$ to $\varphi_j$ if $\varphi_j$ ($2\le j\le r-1$) is a 1-face incident with a (1,1)-edge.\\[3pt]
\textbf{Rule R6.}& $\varphi$ sends $\frac16$ to $\varphi_r$ if $\varphi_r$ is empty and $r\ge 4$.\\[3pt]
\textbf{Rule R7.}& $\varphi$ sends $\frac16$ to $\alpha_2$ if $\alpha$ is incident with a (1,2)-edge and $\alpha_2$ is an empty 0-face.\\[3pt]
\textbf{Rule R8.}& $\varphi$ sends $\frac16$ to $\beta_2$, where $\beta$ is a 2-face of $H_{3-i}$ having exactly one common vertex with 
 $\varphi$ and incident with two (1,2)-edges and $\beta_2$ is an empty 0-face of $H_{3-i}$ adjacent to $\beta$.\\
\end{tabular}}\\[2mm]

\begin{center}
\begin{minipage}[t]{0.32\textwidth}
\begin{center}
\begin{tikzpicture}[scale=0.07]
\draw [thick] (-5,0) -- (45,0);
\draw (20,4) node {$\varphi$};
\draw (10,-4) node {$\alpha$};
\draw (30,-4) node {$\beta$};
\draw (20,21) node {e0-f $\varphi_2$};
\draw [thick,red,->] (20,8) -- (20,15);
\draw [red](24,14) node {$\frac13$};
\draw [thick] (30cm,0cm) arc [start angle=0, end angle=180, radius=10cm];
\draw [thick] (0cm,0cm) arc [start angle=180, end angle=360, radius=10cm];
\draw [thick] (20cm,0cm) arc [start angle=180, end angle=360, radius=10cm];
\filldraw (0,0) circle (7mm);
\filldraw (10,0) circle (7mm);
\filldraw (20,0) circle (7mm);
\filldraw (30,0) circle (7mm);
\filldraw (40,0) circle (7mm);
\filldraw [white] (20,-20) circle (5mm);
\draw (20,-17) node {R4};
\end{tikzpicture}
\end{center}
\end{minipage}
\begin{minipage}[t]{0.32\textwidth}
\begin{center}
\begin{tikzpicture}[scale=0.07]
\draw [thick] (-5,0) -- (45,0);
\draw (20,4) node {$\varphi$};
\draw (10,-4) node {$\alpha$};
\draw (30,-4) node {$\beta$};
\draw (30,10) node {$\varphi_2$};
\draw (25,26) node {e0-f $\varphi_3$};
\draw [thick,red,->] (20,8) -- (25,20);
\draw [red](29,19) node {$\frac13$};
\draw [thick] (30cm,0cm) arc [start angle=0, end angle=180, radius=10cm];
\draw [thick] (40cm,0cm) arc [start angle=0, end angle=180, radius=15cm];
\draw [thick] (0cm,0cm) arc [start angle=180, end angle=360, radius=10cm];
\draw [thick] (20cm,0cm) arc [start angle=180, end angle=360, radius=10cm];
\filldraw (0,0) circle (8mm);
\filldraw (10,0) circle (8mm);
\filldraw (20,0) circle (8mm);
\filldraw (30,0) circle (8mm);
\filldraw (40,0) circle (8mm);
\filldraw [white] (20,-20) circle (5mm);
\draw (20,-17) node {R4};
\end{tikzpicture}
\end{center}
\end{minipage}
\begin{minipage}[t]{0.27\textwidth}
\begin{center}
\begin{tikzpicture}[scale=0.07]
\draw [thick] (-5,0) -- (35,0);
\draw [thick] (30,0) -- (30,-3);
\draw (24,-3) node {1-f};
\draw (20,4) node {$\varphi$};
\draw (10,-4) node {$\alpha$};
\draw (20,21) node {e0-f $\varphi_2$};
\draw [thick,red,->] (20,8) -- (20,15);
\draw [red](24,14) node {$\frac16$};
\draw [thick] (30cm,0cm) arc [start angle=0, end angle=180, radius=10cm];
\draw [thick] (0cm,0cm) arc [start angle=180, end angle=360, radius=10cm];
\filldraw (0,0) circle (8mm);
\filldraw (10,0) circle (8mm);
\filldraw (20,0) circle (8mm);
\filldraw (30,0) circle (8mm);
\filldraw [white] (20,-20) circle (5mm);
\draw (20,-17) node {R4};
\end{tikzpicture}
\end{center}
\end{minipage}
\end{center}
\begin{center}
\begin{minipage}[t]{0.28\textwidth}
\begin{center}
\begin{tikzpicture}[scale=0.062]
\draw [thick] (-5,0) -- (45,0);
\draw [thick] (30,0) -- (30,-3);
\draw (24,-3) node {1-f};
\draw (20,4) node {$\varphi$};
\draw (10,-4) node {$\alpha$};
\draw (30,10) node {$\varphi_2$};
\draw (25,26) node {e0-f $\varphi_3$};
\draw [thick,red,->] (20,8) -- (25,20);
\draw [red](29,19) node {$\frac16$};
\draw [thick] (30cm,0cm) arc [start angle=0, end angle=180, radius=10cm];
\draw [thick] (40cm,0cm) arc [start angle=0, end angle=180, radius=15cm];
\draw [thick] (0cm,0cm) arc [start angle=180, end angle=360, radius=10cm];
\filldraw (0,0) circle (7mm);
\filldraw (10,0) circle (7mm);
\filldraw (20,0) circle (7mm);
\filldraw (30,0) circle (7mm);
\filldraw (40,0) circle (7mm);
\filldraw [white] (20,-20) circle (5mm);
\draw (20,-15) node {R4};
\end{tikzpicture}
\end{center}
\end{minipage}
\begin{minipage}[t]{0.33\textwidth}
\begin{center}
\begin{tikzpicture}[scale=0.062]
\draw [thick] (-5,0) -- (0,0);
\draw [dotted,thick] (0,0) -- (10,0);
\draw [thick] (10,0) -- (30,0);
\draw [dotted,thick] (30,0) -- (40,0);
\draw [thick] (40,0) -- (55,0);
\draw [thick] (20,0) -- (20,-3);
\draw [thick] (40,0) -- (40,-3);
\draw [thick] (50,0) -- (50,-3);
\draw (45,-3) node {1-f};
\draw (20,4) node {$\varphi$};
\draw (43,10) node {$\varphi_j$};
\draw [thick,red,->] (20,8) -- (35,18);
\draw [red](38,16) node {$\frac16$};
\draw [thick] (30cm,0cm) arc [start angle=0, end angle=180, radius=10cm];
\draw [thick] (40cm,0cm) arc [start angle=0, end angle=180, radius=20cm];
\draw [thick] (50cm,0cm) arc [start angle=0, end angle=180, radius=25cm];
\filldraw (0,0) circle (7mm);
\filldraw (10,0) circle (7mm);
\filldraw (20,0) circle (7mm);
\filldraw (30,0) circle (7mm);
\filldraw (40,0) circle (7mm);
\filldraw (50,0) circle (7mm);
\filldraw [white] (20,-20) circle (5mm);
\draw (20,-15) node {R5};
\end{tikzpicture}
\end{center}
\end{minipage}
\begin{minipage}[t]{0.33\textwidth}
\begin{center}
\begin{tikzpicture}[scale=0.062]
\draw [thick] (-5,0) -- (0,0);
\draw [dotted,thick] (0,0) -- (10,0);
\draw [thick] (10,0) -- (40,0);
\draw [dotted,thick] (40,0) -- (50,0);
\draw [thick] (50,0) -- (55,0);
\draw [thick] (20,0) -- (20,-3);
\draw (20,4) node {$\varphi$};
\draw (30,10) node {$\varphi_2$};
\draw (25,36) node {e0-f $\varphi_r$};
\draw [thick,red,->] (20,8) -- (25,30);
\draw [red](29,29) node {$\frac16$};
\draw [thick] (30cm,0cm) arc [start angle=0, end angle=180, radius=10cm];
\draw [thick] (40cm,0cm) arc [start angle=0, end angle=180, radius=15cm];
\draw [thick] (50cm,0cm) arc [start angle=0, end angle=180, radius=25cm];
\filldraw (0,0) circle (7mm);
\filldraw (10,0) circle (7mm);
\filldraw (20,0) circle (7mm);
\filldraw (30,0) circle (7mm);
\filldraw (40,0) circle (7mm);
\filldraw (50,0) circle (7mm);
\filldraw [white] (20,-20) circle (5mm);
\draw (20,-15) node {R6};
\end{tikzpicture}
\end{center}
\end{minipage}
\end{center}
\begin{center}
\begin{minipage}[t]{0.3\textwidth}
\begin{center}
\begin{tikzpicture}[scale=0.07]
\draw [thick] (-5,0) -- (35,0);
\draw [thick] (0,0) -- (0,3);
\draw (5,3) node {1-f};
\draw (20,4) node {$\varphi$};
\draw (10,-4) node {$\alpha$};
\draw (10,-21) node {e0-f $\alpha_2$};
\draw [thick,red,->] (15,4) -- (10,-15);
\draw [red](14,-14) node {$\frac16$};
\draw [thick] (30cm,0cm) arc [start angle=0, end angle=180, radius=10cm];
\draw [thick] (0cm,0cm) arc [start angle=180, end angle=360, radius=10cm];
\filldraw (0,0) circle (7mm);
\filldraw (10,0) circle (7mm);
\filldraw (20,0) circle (7mm);
\filldraw (30,0) circle (7mm);
\filldraw[white] (20,-40) circle (5mm);
\draw (20,-30) node {R7};
\end{tikzpicture}
\end{center}
\end{minipage}
\begin{minipage}[t]{0.36\textwidth}
\begin{center}
\begin{tikzpicture}[scale=0.07]
\draw [thick] (-5,0) -- (55,0);
\draw [thick] (40,0) -- (40,3);
\draw [thick] (50,0) -- (50,3);
\draw (24,-3) node {1-f};
\draw (35,3) node {1-f};
\draw (45,3) node {1-f};
\draw (20,4) node {$\varphi$};
\draw (10,-4) node {$\alpha$};
\draw (40,-4) node {$\beta$};
\draw (40,-21) node {e0-f $\beta_2$};
\draw [thick,red,->] (26,3) to [out=0,in=90] (40,-15);
\draw [red](44,-14) node {$\frac16$};
\draw [thick] (30cm,0cm) arc [start angle=0, end angle=180, radius=10cm];
\draw [thick] (0cm,0cm) arc [start angle=180, end angle=360, radius=10cm];
\draw [thick] (30cm,0cm) arc [start angle=180, end angle=360, radius=10cm];
\filldraw (0,0) circle (7mm);
\filldraw (10,0) circle (7mm);
\filldraw (20,0) circle (7mm);
\filldraw (30,0) circle (7mm);
\filldraw (40,0) circle (7mm);
\filldraw (50,0) circle (7mm);
\filldraw [white] (20,-40) circle (5mm);
\draw (20,-30) node {R8};
\end{tikzpicture}
\end{center}
\end{minipage}
\end{center}

\vspace{-10mm}
\begin{center}
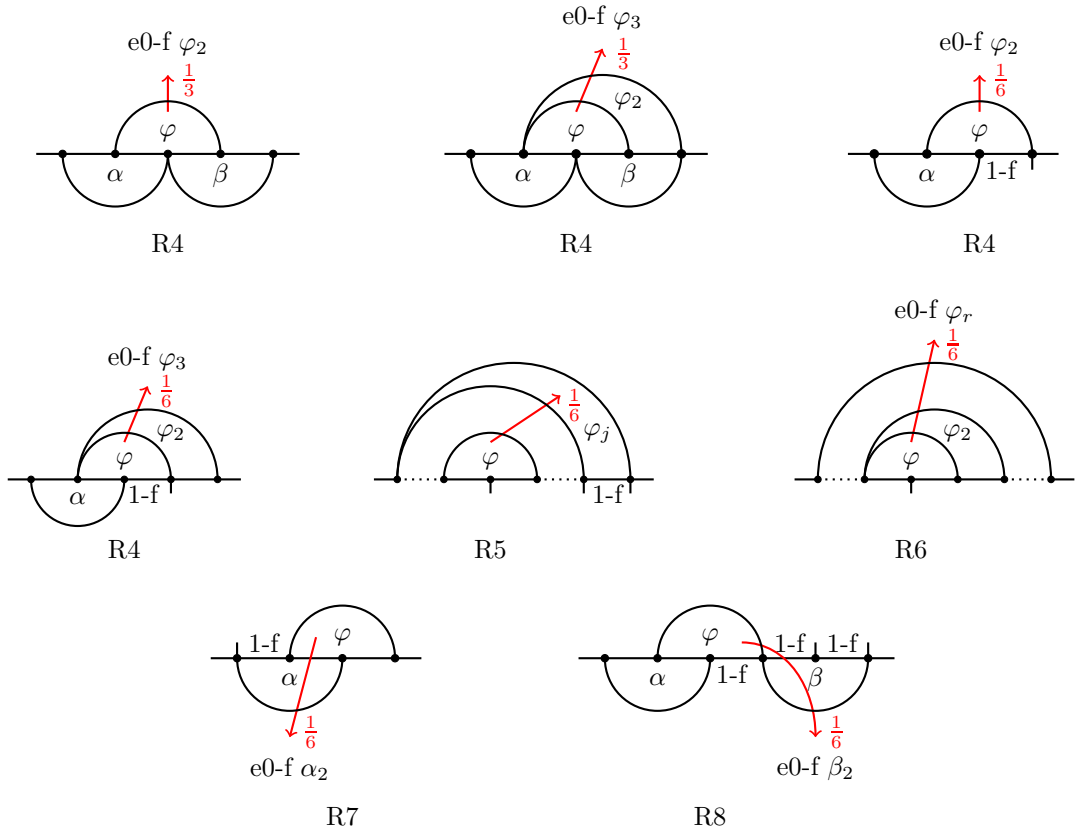
Fig.~2. Redistribution rules R4--8 (1-f is a 1-face and e0-f is an empty 0-face).
\end{center}

\medskip
For an empty face $\varphi$, let $w_2(\varphi)$ be the total weight obtained by $\varphi$ (after second redistribution). Obviously, 
$\sum\limits_{\varphi\in\eF(H)} w_2(\varphi)=|E(C)|=k$ (as non-empty faces do not obtain any weight). In the following, we will show that the weight 
$w_2(\varphi)$ of each (empty) face $\varphi$ does not exceed $\frac23$ which will mean $k=\sum\limits_{\varphi\in\eF(H)} w_2(\varphi)\le \frac23 
\ef(H)$. The maximal planar graph $G$ has exactly $2n-4$ faces. Each of $\ef(H)\ge \frac32k$ empty faces of $H$ is a face of $G$ as well, and each 
of $n-k$ (pairwise non-adjacent) vertices of $G$ not belonging to $C$ (whose removal has created a non-empty face of $H$) is incident with three 
(``private'') faces of $G$. Hence $2n-4=|F(G)|=\ef(H)+3(n-k)\ge \frac32k+3n-3k$ and finally $k\ge \frac23(n+4)$ will follow. 

\medskip
\noindent
\textbf{Weight of a 2-face.}\\[2mm]
Let $\varphi=[v_1,v_2,v_3]$ be a $2$-face of $H_i$ and let $B_{\varphi}=(\varphi,\varphi_2,\dots,\varphi_r)$, $r\ge 2$, be the $\varphi$-branch. 
As already mentioned, $\frac23\le w_1(\varphi)\le 1$. We check that the weight of $\varphi$ exceeding $\frac23$ will be shifted in the second 
redistribution.

\noindent
\textbf{1.} Let $\varphi$ be incident with two (2,2)-edges (note that $w_1(\varphi)=1$). Denote $\alpha=[v_0,v_1,v_2]$ and $\beta=[v_2,v_3,v_4]$ the 
2-faces of $H_{3-i}$ adjacent to $\varphi$. Let $\alpha_2$ and $\beta_2$ be the face of $H_{3-i}$ adjacent to $\alpha$ and $\beta$, respectively. 
Each of the faces $\varphi_2$, $\alpha_2$, and $\beta_2$ is either a 1-face or empty 0-face (by Claim~\ref{empty0}a).

\medskip
\noindent
\textbf{1.1.} Let $\alpha_2$ and $\beta_2$ be 0-faces (possibly $\alpha_2=\beta_2$).

\noindent
\textbf{1.1.1.} If edges $v_0v_1$ and $v_3v_4$ of $C$ do not belong to the rim $R(B_\varphi)$ of $B_\varphi$, then $r=2$, thus $\varphi$ sends 
$\frac13$ to empty 0-face $\varphi_2$ (by R4).

\noindent
\textbf{1.1.2.} If $v_0v_1$ belongs to the rim $R(B_\varphi)$ and $v_3v_4$ does not belong to $R(B_\varphi)$, then $\varphi_2=[v_0,v_1,v_3]$ is 
a 1-face and $\varphi_r$ is empty (by Claim~\ref{empty0}a). Thus $\varphi$ sends weight $\ge\frac16$ to $\varphi_r$ (by R4 or R6) and $\frac16$ 
to $\alpha_2$ (by R7). (Similarly if $v_0v_1$ does not belong to $R(B_\varphi)$ and $v_3v_4$ belongs to $R(B_\varphi)$.)

\noindent
\textbf{1.1.3.} If edges $v_0v_1$ and $v_3v_4$ belong to the rim $R(B_\varphi)$, then both are (1,2)-edges. Thus $\varphi$ sends $\frac16$ to 
$\alpha_2$ and $\frac16$ to $\beta_2$ (by R7).

\medskip
\noindent
\textbf{1.2.} Let $\alpha_2=[v_{-1},v_0,v_2]$ be a 1-face and $\beta_2$ be a 0-face. (Similarly if $\alpha_2$ is a 0-face and $\beta_2$ is a 1-face.)

\noindent
\textbf{1.2.1.} If $v_3v_4$ does not belong to the rim $R(B_\varphi)$, then $r\le 3$ and $\varphi_r$ is empty (by proof of Claim~\ref{empty0}c). 
Thus $\varphi$ sends $\frac13$ to $\varphi_r$ (by R4).

\noindent
\textbf{1.2.2.} If $v_3v_4$ belongs to the rim $R(B_\varphi)$ and $v_0v_1$ does not belong to $R(B_\varphi)$, then $\varphi_2=[v_1,v_3,v_4]$ is 
a 1-face and $\varphi_r$ is empty (by Claim~\ref{empty0}a). Thus $\varphi$ sends weight $\ge\frac16$ to $\varphi_r$ (by R4 or R6) and $\frac16$ 
to $\beta_2$ (by R7).

\noindent
\textbf{1.2.3.} Let edges $v_3v_4$ and $v_0v_1$ belong to the rim $R(B_\varphi)$, then both are (1,2)-edges. If $v_0v_1$ and $v_3v_4$ are 
incident with $\varphi_2$ and $\varphi_3$, then $\{v_0,v_2,v_4\}$ is a non-trivial 3-cut, a contradiction. If $\varphi_2=[v_0,v_1,v_3]$ and 
$\varphi_3=[v_{-1},v_0,v_3]$, then $\{v_{-1},v_2,v_3\}$ is a non-trivial 3-cut, a contradiction as well. Thus $\varphi_2=[v_1,v_3,v_4]$ and 
$\varphi_3=[v_1,v_4,v_5]$.

\noindent
\textbf{1.2.3.1.} If $v_{-1}v_0$ does not belong to the rim $R(B_\varphi)$, then $\varphi_r$ is empty (by Claim~\ref{empty0}b). Thus $\varphi$ 
sends $\frac16$ to $\varphi_r$ (by R6) and $\frac16$ to $\beta_2$ (by R7).

\noindent
\textbf{1.2.3.2.} If $v_{-1}v_0$ belongs to the rim $R(B_\varphi)$, then $v_{-1}v_0$ is a (1,1)-edge. Thus $\varphi$ sends $\frac16$ to $\varphi_j$, 
a 1-face of $B_\varphi$ incident with $v_{-1}v_0$ (by R5) and $\frac16$ to $\beta_2$ (by R7).

\medskip
\noindent
\textbf{1.3.} Let $\alpha_2=[v_{-1},v_0,v_2]$ and $\beta_2=[v_2,v_4,v_5]$ be 1-faces.

\noindent
\textbf{1.3.1.} If $v_3v_4$ does not belong to the rim $R(B_\varphi)$, then $r\le 3$ and $\varphi_r$ is empty (by proof of Claim~\ref{empty0}c). 
Thus $\varphi$ sends $\frac13$ to $\varphi_r$ (by R4). (Similarly if $v_0v_1$ does not belong to $R(B_\varphi)$.)

\noindent
\textbf{1.3.2.} Let edges $v_0v_1$ and $v_3v_4$ belong to the rim $R(B_\varphi)$, then both are (1,2)-edges. If $v_0v_1$ and $v_3v_4$ are 
incident with $\varphi_2$ and $\varphi_3$, then $\{v_0,v_2,v_4\}$ is a non-trivial 3-cut, a contradiction. If $\varphi_2=[v_0,v_1,v_3]$ and 
$\varphi_3=[v_{-1},v_0,v_3]$, then $\{v_{-1},v_2,v_3\}$ is a non-trivial 3-cut, a contradiction as well. (Similarly if $\varphi_2=[v_1,v_3,v_4]$ 
and $\varphi_3=[v_1,v_4,v_5]$.)

\medskip
\noindent
\textbf{2.} Let $\varphi$ be incident with (2,2)-edge $v_1v_2$ and (1,2)-edge $v_2v_3$ (note that $w_1(\varphi)=\frac56$). Denote $\alpha=[v_0,v_1,v_2]$ 
the 2-face of $H_{3-i}$ adjacent to $\varphi$ and let $\alpha_2$ be the face of $H_{3-i}$ adjacent to $\alpha$. Each of the faces $\varphi_2$ and 
$\alpha_2$ is either a 1-face or empty 0-face (by Claim~\ref{empty0}a).

\medskip
\noindent
\textbf{2.1.} Let $\alpha_2$ be 0-face.

\noindent
\textbf{2.1.1.} If $v_0v_1$ does not belong to the rim $R(B\varphi)$, then $\varphi_r$ is empty (by Claim~\ref{empty0}a). Thus $\varphi$ sends 
$\frac16$ to $\varphi_r$ (by R4 or R6).

\noindent
\textbf{2.1.2.} If $v_0v_1$ belongs to the rim $R(B\varphi)$, then $v_0v_1$ is a (1,2)-edge. Thus $\varphi$ sends $\frac16$ to $\alpha_2$ (by R7).

\medskip
\noindent
\textbf{2.2.} Let $\alpha_2$ be a 1-face incident with $v_{-1}v_0$ (i.e. $\alpha_2=[v_{-1},v_0,v_2]$).

\noindent
\textbf{2.2.1.} If $v_3v_4$ does not belong to the rim $R(B_\varphi)$, then $r\le 3$ and $\varphi_r$ is empty (by proof of Claim~\ref{empty0}c). 
Thus $\varphi$ sends $\frac16$ to $\varphi_r$ (by R4).

\noindent
\textbf{2.2.2.} If $v_3v_4$ belongs to the rim $R(B_\varphi)$ and $v_0v_1$ does not belong to $R(B_\varphi)$, then $\varphi_2=[v_1,v_3,v_4]$ is 
a 1-face and $\varphi_r$ is empty (by Claim~\ref{empty0}a). Thus $\varphi$ sends $\frac16$ to $\varphi_r$ (by R4 or R6).

\noindent
\textbf{2.2.3.} Let edges $v_3v_4$ and $v_0v_1$ belong to the rim $R(B_\varphi)$. If $v_{-1}v_0$ does not belong to $R(B_\varphi)$, then $\varphi_r$ 
is empty (by Claim~\ref{empty0}b). Thus $\varphi$ sends $\frac16$ to $\varphi_r$ (by R6). Otherwise $v_{-1}v_0$ belongs to  $R(B_\varphi)$, thus it 
is a (1,1)-edge incident with a 1-face $\varphi_j$ of $B_\varphi$. Hence $\varphi$ sends $\frac16$ to $\varphi_j$ (by R5).

\medskip
\noindent
\textbf{2.3.} Let $\alpha_2$ be a 1-face incident with $v_2v_3$ (i.e. $\alpha_2=[v_0,v_2,v_3]$). Since $v_0v_3\in E(H_{3-i})$, $\varphi_2$ cannot be 
the 1-face $[v_0,v_1,v_3]$ in $H_i$.

\noindent
\textbf{2.3.1.} If $v_3v_4$ does not belong to the rim $R(B_\varphi)$, then $r=2$, thus $\varphi$ sends $\frac16$ to $\varphi_2$ (by R4).

\noindent
\textbf{2.3.2.} If $v_3v_4$ belongs to the rim $R(B_\varphi)$, then $r\ge 3$ and $\varphi_2=[v_1,v_3,v_4]$.

\noindent
\textbf{2.3.2.1.} If $v_3v_4$ is incident with a 1-face of $H_{3-i}$ (i.e., $v_3v_4$ is a (1,1)-edge), then $\varphi$ sends $\frac16$ to $\varphi_2$ 
(by~R5).

\noindent
\textbf{2.3.2.2.} Let $v_3v_4$ be incident with a 2-face $\beta$ of $H_{3-i}$ (necessarily, $\beta=[v_3,v_4,v_5]$). 
If $r=3$, then $\varphi_3$ is empty (by Claim~\ref{empty0}a), thus $\varphi$ sends $\frac16$ to $\varphi_3$ (by R4). 
If $r=4$, then $\varphi_3=[v_1,v_4,v_5]$ (as $\{v_0,v_3,v_4\}$ is a non-trivial 3-cut if $\varphi_3=[v_0,v_1,v_4]$) and $\varphi_4$ is empty 
(by Claim~\ref{empty0}a), thus $\varphi$ sends $\frac16$ to $\varphi_4$ (by R6). 
Finally, let $r\ge 5$. Necessarily $\varphi_3=[v_1,v_4,v_5]$ (as for $r=4$) and $\varphi_4=[v_1,v_5,v_6]$ (as $\{v_0,v_3,v_5\}$ is a non-trivial 
3-cut if $\varphi_4=[v_0,v_1,v_5]$) are 1-faces of $B_\varphi$. If $v_5v_6$ is a (1,1)-edge, then $\varphi$ sends $\frac16$ to $\varphi_4$ 
(by R5). Otherwise $v_5v_6$ is a (1,2)-edge, thus it does not belong to $\beta$-branch (in $H_{3-i}$) and therefore $\beta_2$ is a 0-face, which is, 
moreover, empty (as the cycle obtained from $C$ by replacing the path $(v_0,\dots,v_5)$ by the path $(v_0,v_2,v_1,v_4,v_3,v_5)$ is a longest good cycle 
of $G$ and contains the edge $v_3v_5$ incident with $\beta_2$ (Claim~\ref{empty12})). Hence $\varphi$ sends $\frac16$ to $\beta_2$ (by R8). 

\medskip
\noindent
\textbf{Weight of a 1-face.}\\[2mm]
To estimate the weight of a 1-face, we use the following simple observation:

\begin{claim}\label{1-face}
Each 1-face of $H$ belongs to at most one branch.
\end{claim}

\noindent
Let $\psi$ be a $1$-face incident with an edge $e$ of $C$. If $e$ is a (1,2)-edge, then $\psi$ obtains weight $\frac23$ from $e$ (by R2) only. 
Otherwise $e$ is a (1,1)-edge, thus $\psi$ obtains $\frac12$ from $e$ (by R1). Furthermore, in this case, $\psi$ can get $\frac16$ from 
a 2-face $\varphi$ (by R5) if $\psi$ belongs to the $\varphi$-branch. Hence $w_2(\psi)\le\frac23$.

\medskip
\noindent
\textbf{Weight of an empty 0-face.}\\[2mm]
Each empty 0-face $\omega$ belongs to at most two branches (in Case 1). Let $\varphi$ be a $2$-face of $H_i$ with the $\varphi$-branch 
$B_\varphi=(\varphi,\varphi_2,\dots,\varphi_r)$ such that $\varphi_r=\omega$, and let $e$ be the edge incident with $\varphi_r$ and $\varphi_{r-1}$ 
(where $\varphi_{r-1}=\varphi$ for $r=2$).

If $\varphi$ is adjacent to two 2-faces, then $\omega$ gets through $e$ the weight $\frac13$ (by R4) for $r\le 3$ or the weight $\frac16$ (by R6) 
for $r\ge 4$. If $\varphi$ is adjacent to one 2-face, then $\omega$ gets through $e$ the weight $\frac16$ (by R4) and additionally $\frac16$ (by R7) 
for $r=2$ or the weight at most $\frac16$ (by R4) for $r=3$ or the weight $\frac16$ (by R6) for $r\ge 4$. Finally, if $\varphi$ is adjacent to no 
2-face, then $\omega$ gets through $e$ the weight $\frac16$ (by R6) for $r\ge 4$ or the weight at most $2\times\frac16$ (by R8) for $r\le 3$.

\medskip
We showed that $w_2(\varphi)\le\frac23$ for each empty face $\varphi$ and completed the Case~1. 
Thus, we can assume that in $H_i$ are only empty faces and among them, at most one face is a 0-face. To complete the proof, we have to show that there 
are some empty faces in $H_{3-i}$ as well.

\medskip
\noindent
\textbf{CASE 2.} Let $H_i$ contain no 0-face or exactly one 0-face which is additionally empty.\\[2mm]
Obviously, if $H_i$ contains no 0-face, then it contains two 2-faces $\alpha_1$ and $\alpha_2$ (since $T_i$ is a path and 2-faces of $H_i$ are leaves 
of $T_i$). Note that, (only) in this case, the branches in $H_i$ are not defined.

Remember that $H=G[V(C)]$ has $k\ge 7$ vertices (as otherwise $G$ with at most $k+2\le 8$ vertices is Hamiltonian). If $H_i$ contains exactly one 
0-face, then it contains three 2-faces $\alpha_1$, $\alpha_2$ and $\alpha_3$ (since $T_i$ is a subdivision of $K_{1,3}$ and 2-faces of $H_i$ are leaves 
of $T_i$). We assume that $H_{3-i}$ contains at least two 0-faces as otherwise all but at most one faces of $H_{3-i}$ are empty and $G$ has 
$n\le |V(H)|+1=k+1$ vertices and Theorem~\ref{main} follows immediately (with $n\ge 11$).

\medskip
\newpage
\noindent
\textbf{Distribution of points.}\\[2mm]
To estimate the number of empty 0- and 1-faces in $H_{3-i}$, each 2-face $\alpha_j$ of $H_i$ ($j\in\{1,2\}$ if $H_i$ contains no 0-face and 
$j\in\{1,2,3\}$ if $H_i$ contains one 0-face, respectively) will distribute 1 or 2 points to faces of $H_{3-i}$. Let $\alpha_j$ be adjacent to 
the faces $\varphi$ and $\psi$ of $H_{3-i}$.

\medskip
\noindent
\hspace*{-7pt}
{\setlength{\tabcolsep}{3pt}
\begin{tabular}{lp{14cm}}
\textbf{Rule P1.}& If $\varphi$ and $\psi$ are 2-faces of $H_{3-i}$ with branches $B_{\varphi}=(\varphi,\varphi_2,\dots,\varphi_r)$ and $B_{\psi}=
(\psi,\psi_2,\dots,\psi_t)$, then $\varphi_r$ and $\psi_t$ will each receive 1 point (or 2 points if $\varphi_r=\psi_t$) from $\alpha_j$.\\[3pt]
\textbf{Rule P2.}& If $\varphi$ and $\psi$ are 1-faces of $H_{3-i}$, then $\varphi$ and $\psi$ will each receive 1 point from $\alpha_j$.\\[3pt]
\textbf{Rule P3.}& If $\varphi$ is a 2-faces of $H_{3-i}$ with $\varphi$-branch $B_{\varphi}=(\varphi,\varphi_2,\dots,\varphi_r)$ and $\psi$ is 
a 1-face of $H_{3-i}$ not belonging to $B_{\varphi}$, then $\varphi_r$ and $\psi$ will each receive 1 point from $\alpha_j$.\\[3pt]
\textbf{Rule P4.}& If $\varphi$ is a 2-faces of $H_{3-i}$ with $\varphi$-branch $B_{\varphi}=(\varphi,\varphi_2,\dots,\varphi_r)$ and $\psi$ is 
a 1-face of $H_{3-i}$ belonging to $B_{\varphi}$, then only $\psi$ will receive 1 point from $\alpha_j$.\\
\end{tabular}

\medskip
For a face $\varphi$ of $H_{3-i}$, let $p(\varphi)$ be the total number of points carried by $\varphi$ (in the distribution of points).

\begin{claim}\label{f1+2f0 ge p}
$f_1(H_{3-i})+2f_0(H_{3-i})\ge \sum\limits_{\varphi\in\eF(H_{3-i})}p(\varphi)$.
\end{claim}

\noindent
\textit{Proof.}
We have to prove that each 1-face of $H_{3-i}$ gets at most 1 point and that each 0-face of $H_{3-i}$ gets points only if it is empty and it gets 
at most 2 points. Consequently, Claim~\ref{f1+2f0 ge p} follows by simple counting.

Let $\beta$ be a 1-face of $H_{3-i}$. Since $\beta$ can only get points if it is adjacent to some $\alpha_j$ and there can only be one such face 
then $p(\beta)\le 1$.

Let $\beta$ be a 0-face of $H_{3-i}$. Since $\beta$ can only get points if it belongs to a branch and it belongs to at most two branches (as there are 
at least two 0-faces in $H_{3-i}$), then $p(\beta)\le 2$. Assume first that $\beta$ gets a point by P1. Then there is $\alpha_j$ incident with two 
$(2,2)$-edges and adjacent 2-faces $\varphi$ and $\psi$ of $H_{3-i}$. Let $B_{\varphi}=(\varphi,\varphi_2,\dots,\varphi_r)$ with $\varphi_r=\beta$ be 
the branch which ends in $\beta$. By Claim~\ref{empty0}a, $\varphi_r=\beta$ is an empty 0-face.

Thus, assume that $\beta$ gets a point by P3. Then there is $\alpha_j$ incident with a $(1,2)$-edge with adjacent 1-face $\psi$ in $H_{3-i}$  and 
a $(2,2)$-edge with adjacent 2-face $\varphi$ such that $\psi$ does not belong to the branch $B_\varphi=(\varphi,\varphi_2,\dots,\varphi_r)$ with 
$\varphi_r=\beta$. Since the common edge of $\alpha_j$ and $\psi$ does not belong to the rim $R(B_\varphi)$, again by Claim~\ref{empty0}a, 
$\varphi_r=\beta$ is an empty 0-face.
\hfill $\Box$

\medskip
\begin{claim}\label{f1+2f0}
$f_1(H_{3-i})+2f_0(H_{3-i})\ge 4$.
\end{claim}

\noindent
\textit{Proof.}
If $\sum\limits_{\varphi\in\eF(H_{3-i})}p(\varphi)\ge 4$, then $f_1(H_{3-i})+2f_0(H_{3-i})\ge 4$ (by Claim~\ref{f1+2f0 ge p}). 
Assume $\sum\limits_{\varphi\in\eF(H_{3-i})}p(\varphi)\le 3$.

\medskip
\noindent
\textbf{1.} Let $H_i$ contains exactly one 0-face. As there are three 2-faces $\alpha_1, \alpha_2, \alpha_3$ in $H_i$
(note, that $T_i$ is a subdivided 3-star in this case), then 
$\sum\limits_{\varphi\in\eF(H_{3-i})}p(\varphi)=3$. Furthermore, only P4 was applied to each $\alpha_j$ ($j\in\{1,2,3\}$) hence there are three 
1-faces with 1 point and they belong to three different branches.

Since $|V(H)|=k\ge 7$, there is $j\in\{1,2,3\}$ such that $\alpha_j$ is adjacent to a 1-face $\delta$ of $H_i$. Let $\varphi$ be the adjacent 2-face 
of $\alpha_j$ in $H_{3-i}$ and $B_\varphi=(\varphi,\varphi_2,\dots,\varphi_r)$ be its branch.

\medskip
\noindent
\textbf{1.1.} If $r\ge 4$, then $\varphi_2$ and $\varphi_3$ are 1-faces of the same branch. Thus, at most one among $\varphi_2$ and $\varphi_3$ has 
a point and $f_1(H_{3-i})\ge 4$.

\medskip
\noindent
\textbf{1.2.} If $r=3$, then $\delta$ and $\varphi$ are not adjacent (i.e. $\delta\neq\varphi_2$, since $H$ has no multiple edges) and $\varphi_3$ is 
an empty 0-face (by Claim~\ref{empty0}b), hence $f_1(H_{3-i})+f_0(H_{3-i})\ge 4$.

\medskip
\newpage
\noindent
\textbf{2.} Let $H_i$ contains no 0-face. Since $\sum\limits_{\varphi\in\eF(H_{3-i})}p(\varphi)\le 3$, there is $j\in\{1,2\}$ such that 
P4 was applied to $\alpha_j$. Let $\delta$ be the 1-face of $H_i$ adjacent to $\alpha_j$ (since $|V(H)|=k\ge 7$), let $\varphi$ and $\psi$ be the 2-face 
and 1-face of $H_{3-i}$ adjacent with $\alpha_j$, respectively, and let $B_\varphi=(\varphi,\varphi_2,\dots,\varphi_r)$ be the branch of $\varphi$.
We may assume $\alpha_j=[v_1,v_2,v_3]$ and $\varphi=[v_2,v_3,v_4]$.

\medskip
\noindent
\textbf{2.1.} Let $r\le 4$.

\noindent
\textbf{2.1.1} If $\delta=[v_0,v_1,v_3]$, then $v_0v_1$ does not belong to the rim $R(B_{\varphi})$ (otherwise $\varphi_2=[v_1,v_2,v_4]$, 
$\varphi=[v_0,v_1,v_4]$ and ${v_0,v_3,v_4}$ is a non-trivial 3-cut, a contradiction) and $\varphi_r$ is an empty 0-face (by Claim~\ref{empty0}b). 
By P1--4, there is a face in $H_{3-i}$ other than $\psi$ and $\varphi_r$ with a point, thus $f_1(H_{3-i})+2f_0(H_{3-i})\ge 4$.

\noindent
\textbf{2.1.2} If $\delta=[v_1,v_3,v_4]$, then $\varphi_2=[v_2,v_4,v_5]$ (since $v_1v_4\in E(H_i)$), $\psi=\varphi_3=[v_1,v_2,v_5]$, and 
$\{v_1,v_4,v_5\}$ is a non-trivial 3-cut, a contradiction.

\medskip
\noindent
\textbf{2.2.} Let $r=5$. There are three 1-faces (in fact $\varphi_2$, $\varphi_3$, and $\varphi_4$) all belonging to the same branch $B_\varphi$. 
We may assume that P4 was applied to $\alpha_j$ and P2 was applied to $\alpha_{3-j}$, and all three 1-faces are adjacent to $\alpha_1$ or $\alpha_2$ 
(since otherwise there is another 1-face or empty 0-face and Claim~\ref{f1+2f0} follows).

\noindent
\textbf{2.2.1.} If $\alpha_{3-j}=[v_{-1},v_0,v_1]$, then rim $R(B_\varphi)=(v_{-1},\dots,v_4)$, thus $\varphi_2=[v_1,v_2,v_4]$ and 
$\delta=[v_1,v_3,v_4]$, a contradiction to the simplicity of $H$.

\noindent
\textbf{2.2.2.} If $\alpha_{3-j}=[v_4,v_5,v_6]$ and $\delta=[v_0,v_1,v_3]$, then rim $R(B_\varphi)=(v_1,\dots,v_6)$ and $\varphi_5$ is an empty 0-face 
(by Claim~\ref{empty0}b), thus $f_1(H_{3-i})+f_0(H_{3-i})\ge 4$.

\noindent
\textbf{2.2.3.} If $\alpha_{3-j}=[v_4,v_5,v_6]$ and $\delta=[v_1,v_3,v_4]$, then rim $R(B_\varphi)=(v_1,\dots,v_6)$. Hence $v_1v_6\in E(H_{3-i})$ 
and consequently $\{v_1,v_4,v_6\}$ is a non-trivial 3-cut, a contradiction.

\medskip
\noindent
\textbf{2.3.} If $r\ge 6$, then there are at least four 1-faces in $B_\varphi$, thus $f_1(H_{3-i})\ge 4$.
\hfill $\Box$

\medskip
Remember that each $j$-face of $H_{3-i}$ is incident with $j$ (``private'') edges of $C$, hence $2f_2(H_{3-i})+f_1(H_{3-i})=k$. As each of the $k-2$ 
triangular faces of $H_i$ is empty, all non-empty faces of $H$ belong to $H_{3-i}$ and their number is $(k-2)-f_2(H_{3-i})-f_1(H_{3-i})-f_0(H_{3-i})=
(k-2)-\frac12(k-f_1(H_{3-i}))-f_1(H_{3-i})-f_0(H_{3-i})=\frac{k}{2}-2-\frac12(f_1(H_{3-i})+2f_0(H_{3-i}))\le \frac{k}{2}-4$ (by Claim~\ref{f1+2f0}). 
Finally, at most $\frac{k}{2}-4$ vertices of $G$ lie outside the cycle $C$ (and exactly $k$ vertices on $C$), hence $n\le k+(\frac{k}{2}-4)$ and $k\ge
\frac23(n+4)$ follows, which completes the proof of Theorem~\ref{main}.

\bibliographystyle{abbrvurl}
\bibliography{FHMS}

\end{document}